# Summation Formulas, from Poisson and Voronoi to the Present

Stephen D. Miller[1]
Wilfried Schmid[2]

In Honor of Jacques Carmona

## 1 Introduction

Summation formulas have played a very important role in analysis and number theory, dating back to the Poisson summation formula. The modern formulation of Poisson summation asserts the equality

$$(1.1) \qquad \sum_{n\in\mathbb{Z}} f(n) \;=\; \sum_{n\in\mathbb{Z}} \widehat{f}(n) \quad \left(\widehat{f}(t) \;=\; \int_{\mathbb{R}} f(x)\, e^{-2\pi i x t}\, dx\right),$$

valid (at least) for all Schwartz functions $f$. Let us take a brief historical detour to the beginning of the 20th century, before the notion of Schwartz function had been introduced. The custom then was to state (1.1) for more general functions $f$, such as functions of bounded variation, but supported on a finite interval, and usually in terms of the cosine:

$$(1.2) \qquad {\sum_{a\leq n\leq b}}' f(n) \;=\; \int_a^b f(x)\, dx \;+\; 2\sum_{n=1}^\infty \int_a^b f(x)\,\cos(2\pi n x)\, dx\,;$$

the notation $\sum'$ signifies that at points $n$ where $f$ has a discontinuity – including the endpoints $a$, $b$ – the term $f(n)$ is to be interpreted as the average of the left and right limits of $f(x)$. Indeed, the general case of (1.2) can be reduced to the special case of $a=0$, $b=1$, which amounts to the statement that the Fourier series of a periodic function of bounded variation converges pointwise, to the average of its left and right-hand limits.

To Voronoi and his contemporaries, formulas like (1.2) appeared tailor-made to evaluate certain finite sums. That seemed significant because several important questions in number theory involve finite sums of arithmetic quantities.

---

[1] Supported by NSF grant DMS-0122799
[2] Supported in part by NSF grant DMS-0070714

In this connection, let us mention two classical examples, *Dirichlet's divisor problem* and the *Gauss' circle problem*. The former estimates the size of $d(n)$, the number of positive divisors of an integer $n$. Dirichlet proved

$$(1.3) \qquad D(X) =_{\text{def}} \sum_{n=1}^{X} d(n) \;-\; X \log X - (2\gamma - 1)X \;=\; O(X^{1/2}),$$

where $\gamma$ is Euler's constant $\approx .57721566$. Gauss' circle problem concerns the average size of

$$(1.4) \qquad r_2(n) =_{\text{def}} \#\{\,(x,y) \in \mathbb{Z}^2 \mid x^2 + y^2 = n\,\},$$

for which Gauss gave the estimate

$$(1.5) \qquad \Delta(X) =_{\text{def}} \left| \sum_{n=1}^{X} r_2(n) - \pi X \right| \;=\; O(X^{1/2}).$$

Each problem has a geometric interpretation, with $D(X)$ counting lattice points in the region $\{\,x, y > 0,\ xy \leq X\}$, and $\Delta(X)$ lattice points in the disc $\{\,x^2 + y^2 \leq X\}$. These two bounds are related, as we shall see, and come from fairly elementary considerations.

In the series of papers [27–29] Voronoi developed geometric and analytic methods to improve both Dirichlet's and Gauss' bound. Most importantly in retrospect, he generalized (1.2) by allowing weighted sums, at the expense of introducing more general integral operations on $f$ than the Fourier transform. His formula for the divisor function reads as follows:

$$(1.6) \qquad \begin{aligned} \sideset{}{'}\sum_{a \leq n \leq b} d(n)\, f(n) \;=&\; \int_a^b f(x)\,(\log x + 2\gamma)\, dx \;+\; \\ &+\; \sum_{n=1}^{\infty} d(n) \int_a^b f(x) \left(4 K_0(4\pi\sqrt{nx}) - 2\pi Y_0(4\pi\sqrt{nx})\right) dx\,; \end{aligned}$$

here $f$ can be any piecewise continuous, piecewise monotone function, and $Y_0$, $K_0$ denote Bessel functions; for details see [12]. Voronoi also asserted a formula for $r_2(n)$, which was later rigorously proved by Hardy [7] and Sierpiński [25]:

$$(1.7) \qquad \sideset{}{'}\sum_{a \leq n \leq b} r_2(n)\, f(n) \;=\; \sum_{n=0}^{\infty} r_2(n) \int_a^b f(x)\, \pi\, J_0(2\pi\sqrt{nx})\, dx,$$

under the same hypotheses on the test function $f$, and with $J_0$ referring to the $J$-Bessel function.

These formulas have been used to improve the error estimates in (1.3) and (1.5) to $O(X^{1/3})$ and beyond, and also in the opposite direction, to show that $D(X)$ and $\Delta(X)$ cannot be bounded by a multiple of $X^{1/4}$. Hardy and Landau made that observation [8], which led them to conjecture that $O(X^{1/4+\epsilon})$ is the best possible estimate for both the Dirichlet and Gauss problems. We shall return to this topic in section 6.

Voronoi saw (1.2) and (1.6–1.7) as fitting into a general pattern: he postulated the existence of analogous formulas for sums $\sum'_{a \leq n \leq b} a_n f(n)$, corresponding to any "arithmetic" sequence of coefficients $a_n$, $1 \leq n < \infty$. The wide latitude in the choice of test functions $f$, including characteristic functions of intervals, seemed like an important feature of these summation formulas. The current point of view is different. Characteristic functions and other discontinuous test functions can readily be approximated by smooth functions, and little information about the left hand side of (1.2) or (1.6–1.7) is lost by doing so. On the other hand, the integral transforms of the smooth approximations decay rapidly, which makes the right hand side of the summation formulas much easier to analyze. We shall demonstrate this by example in section 6, where we show how to obtain the estimate $O(X^{1/3+\epsilon})$ in Gauss' circle problem directly from Poisson summation.

In an another direction, Poisson summation gave birth to modular identities, such as the $\theta$-identity

$$(1.8) \qquad \theta(t) \;=_{\text{def}}\; \sum_{n \in \mathbb{Z}} e^{-\pi n^2 t} \;=\; \frac{1}{\sqrt{t}}\, \theta(1/t) \qquad (t > 0),$$

which reflects the fact that $e^{-\pi x^2}$ is its own Fourier transform. Other modular forms, and also $L$-functions, can be created from Poisson's and Voronoi's formulas. Indeed, Riemann developed the theory of his zeta function starting from (1.8), to which he applied the Mellin transform. There are implications in both directions: nowadays the coefficients of $L$-functions are regarded as the natural class of coefficients $(a_n)_{n \geq 1}$ for which Voronoi summation formulas can be proved; in fact, the summation formulas are deduced from properties of the $L$-function in question.

The mechanism of deriving Voronoi summation from $L$-functions is well understood for modular forms and Maass forms on the upper half plane, but runs into difficulties for more general types of automorphic forms. Recently we were able to overcome these difficulties by representation-theoretic arguments: our new result is a Voronoi formula for automorphic forms on $GL(3)$, which we shall state in section 5. It seems clear that our methods extend to automorphic forms on $GL(n)$, for any $n$. Interestingly, the first application of our $GL(3)$ formula involves Maass forms on $GL(2)$, via the symmetric-square lift. We shall discuss the relevant result, due to Sarnak and Watson, following the statement of our formula.

During the genesis of this note, we learned much about the history of the Voronoi summation formula in conversation with several colleagues. Our thanks go especially to Noam Elkies, Martin Huxley, Henryk Iwaniec, Philippe Michel, Peter Sarnak, and Trevor Wooley.

## 2   Poisson summation and the $\zeta$-function

André Weil has written a fascinating historical account of the early history of the Riemann zeta function [31, 32]. At the formal level, it turns out, the connection between the Poisson summation formula and the function $\zeta(s) = \sum_{n=1}^{\infty} n^{-s}$ was recognized even before Riemann. For the statement of Riemann's celebrated result, let us recall the definition of the related function

$$\xi(s) \;=\; \pi^{-s/2}\, \Gamma(s/2)\, \zeta(s) \tag{2.1}$$

and the notion of *order* of an entire function $f(s)$ – the least positive number $\mu$ such that $|f(s)| = O(e^{|s|^{\mu+\epsilon}})$, for every $\epsilon > 0$.

**2.2 Theorem.** (Riemann [22]) *The function $\xi(s)$ has a meromorphic continuation to $\mathbb{C}$, which obeys the functional equation $\xi(s) = \xi(1-s)$. It has simple poles at $s = 0$ and $s = 1$, and the product $s(s-1)\xi(s)$ is entire, of order 1.*

In this section we shall explain the "equivalence" between theorem 2.2 and the Poisson summation formula for Schwartz functions. Both results stand on their own, of course; the point is that each can be derived from the other. In one direction, Riemann proved theorem 2.2 using Poisson summation. More precisely, he applied the Mellin transform to the $\theta$-relation (1.8), which follows from the Poisson summation formula (1.1), with $f(x) = e^{-\pi t x^2}$; details can be found in [3, 5, 22], for example.

The opposite implication is classical, too, but less well-known. We sketch the argument because it can serve as model for the derivation of Voronoi summation from properties of the corresponding $L$-function. First we observe that it suffices to deduce (1.1) from theorem 2.2 for even functions only, since both sides of the equation vanish when $f$ is an odd function. In other words, we must establish the identity

$$\sum_{n=1}^{\infty} f(n) \;=\; \int_0^{\infty} f(t)\, dt \;-\; \frac{f(0)}{2} \;+\; \sum_{n=1}^{\infty} F(n)\,,$$
$$\text{where} \quad F(y) \;=\; \int_{\mathbb{R}} f(x) \cos(2\pi x y)\, dx\,, \tag{2.3}$$

for every even Schwartz function $f$. In that case the Mellin transform of $f$,

$$Mf(s) \;=\; \int_0^{\infty} f(t)\, t^{s-1}\, dt\,, \tag{2.4}$$

is holomorphic in the region $\operatorname{Re} s > -2$, except for a simple pole at $s = 0$, with residue $f(0)$. The Mellin inversion formula

$$f(x) \;=\; \frac{1}{2\pi i} \int_{\operatorname{Re} s = \sigma} Mf(s)\, |x|^{-s}\, ds \qquad (\sigma > 0) \tag{2.5}$$

allows us to express the left hand side of (2.3) as

$$(2.6) \qquad \sum_{n=1}^{\infty} f(n) \;=\; \frac{1}{2\pi i} \int_{\operatorname{Re} s=\sigma} \zeta(s)\, Mf(s)\, ds \qquad (\sigma > 1).$$

When we shift the contour of integration in this formula from $\operatorname{Re} s = \sigma > 1$ to $\operatorname{Re} s = -1$, we pick up a pole at $s=1$ with residue $Mf(1)$, and a pole at $s=0$ with residue $\zeta(0)f(0) = -1/2\, f(0)$. These are the first two terms on the right hand side of (2.3). Note that the contour shift is justified: $\zeta(s)$ grows slowly and $Mf(s)$ decays rapidly as $|\operatorname{Im} s| \to \infty$. We now change variables from $s$ to $1-s$ and apply the functional equation:

$$(2.7) \qquad \begin{aligned} \sum_{n=1}^{\infty} f(n) - Mf(1) + \frac{f(0)}{2} &= \frac{1}{2\pi i} \int_{\operatorname{Re} s=-1} \zeta(s)\, Mf(s)\, ds \\ &= \frac{1}{2\pi i} \int_{\operatorname{Re} s=2} \frac{\pi^{-s/2}\, \Gamma(\frac{s}{2})}{\pi^{(s-1)/2}\, \Gamma(\frac{1-s}{2})} \zeta(s)\, Mf(1-s)\, ds. \end{aligned}$$

But

$$(2.8) \qquad MF(s) \;=\; \frac{\pi^{-s/2}\, \Gamma(\frac{s}{2})}{\pi^{(s-1)/2}\, \Gamma(\frac{1-s}{2})}\, Mf(1-s)$$

is the Mellin transform of the function $F$ in the identity (2.3). This follows from the classical identity

$$(2.9) \qquad \int_{\mathbb{R}} \cos(2\pi x)\, |x|^{s-1}\, dx \;=\; \frac{\pi^{-s/2}\, \Gamma(\frac{s}{2})}{\pi^{(s-1)/2}\, \Gamma(\frac{1-s}{2})} \qquad (0 < \operatorname{Re} s < 1),$$

initially in the region $\{0 < \operatorname{Re} s < 1\}$, but then for all $s \in \mathbb{C}$ by meromorphic continuation. The even Schwartz function $F$ can play the part of $f$ in (2.6). At this point, (2.7–2.8) imply (2.3).

## 3    $L$-functions and Voronoi summation

The two Dirichlet series with coefficients $d(n)$ and $r_2(n)$, respectively, also satisfy functional equations. In the case of $d(n)$,

$$(3.1) \qquad \sum_{n=1}^{\infty} d(n)\, n^{-s} \;=\; \sum_{m,n \geq 1} (m\,n)^{-s} \;=\; \zeta^2(s)$$

inherits its functional equation and analytic properties from the Riemann zeta function. The generating series for $r_2(n)$,

$$(3.2) \qquad \zeta_{\mathbb{Q}(\sqrt{-1})}(s) \;=\; \frac{1}{4} \sum_{n=1}^{\infty} r_2(n)\, n^{-s},$$

is known as the Dedekind zeta function of the number field $\mathbb{Q}(\sqrt{-1})$, because it counts the number of ideals of norm $n$. It factors as the product of two Dirichlet series:

$$\zeta_{\mathbb{Q}(\sqrt{-1})}(s) \;=\; \zeta(s)\, L(s, \chi_4)\,, \tag{3.3}$$

where

$$L(s, \chi_4) \;=\; \sum_{k=0}^{\infty} \frac{(-1)^k}{(2k+1)^s} \tag{3.4}$$

is the first example of a Dirichlet $L$-function. It is entire, of order 1, and satisfies the functional equation

$$L(s, \chi_4) \;=\; 2^{1-2s}\, \frac{\pi^{-(2-s)/2}\, \Gamma(\frac{2-s}{2})}{\pi^{-(s+1)/2}\, \Gamma(\frac{s+1}{2})}\, L(1-s, \chi_4)\,. \tag{3.5}$$

The factorization (3.3) is a result from class field theory, which reflects the fact that an odd prime can be expressed as the sum of two squares if and only if it is congruent to 1 modulo 4. The factorization is also equivalent to the identity

$$\begin{aligned} r_2(n) \;&=\; 4\, \#\{\, d \in \mathbb{N} \mid d|n,\ d \equiv 1 \bmod 4\,\} \\ &\quad -\, 4\, \#\{\, d \in \mathbb{N} \mid d|n,\ d \equiv -1 \bmod 4\,\}\,, \end{aligned} \tag{3.6}$$

which implies the bound

$$r_2(n) \;\leq\; 4\, d(n) \;=\; O(n^\epsilon), \quad \text{for any } \epsilon > 0\,; \tag{3.7}$$

this uses a standard estimate for $d(n)$, which can be found in [21, p. 186], for example.

The Voronoi summation formulas (1.6) and (1.7) can be derived from the functional equations of the Dirichlet series (3.1) and (3.2), much in the same way as Poisson summation was derived from theorem 2.2. To match the $\Gamma$-factors in the functional equations to the integral kernels on the right hand sides of (1.6–1.7), one uses the identities

$$\begin{aligned} \int_0^\infty \left(4 K_0(\sqrt{x}) - 2\pi Y_0(\sqrt{x})\right) x^{s-1} dx \;&= \\ &=\; 2^{4s} \pi^{2s} \left( \frac{\pi^{-s/2}\, \Gamma(\frac{s}{2})}{\pi^{(s-1)/2}\, \Gamma(\frac{1-s}{2})} \right)^2, \\ \int_0^\infty J_0(\sqrt{x}) x^{s-1} dx \;&=\; 4^s\, \frac{\Gamma(s)}{\Gamma(1-s)}\,, \end{aligned} \tag{3.8}$$

which are valid for $0 < \operatorname{Re} s < 3/4$. In effect, these identities take the place of (2.9) in the argument in the previous section.

Modular forms are a rich source of $L$-functions, and each such $L$-function satisfies a functional equation from which a Voronoi formula can be deduced.

To simplify the discussion, we shall consider the case of "full level", i.e., modular forms invariant under all of $SL(2,\mathbb{Z})$. Recall the definition of modular form of weight $k$: a holomorphic function $F : \{\operatorname{Im} z > 0\} \to \mathbb{C}$ such that

$$
\text{(3.9)} \quad F(z) = (cz+d)^{-k} F\left(\frac{az+b}{cz+d}\right) \text{ if } \begin{pmatrix} a & b \\ c & d \end{pmatrix} \in SL(2,\mathbb{Z})
$$
$$
\text{and } |F(z)| = O((\operatorname{Im} z)^N) \text{ as } \operatorname{Im} z \to \infty,
$$

for some $N > 0$. The invariance condition implies invariance under $z \mapsto z+1$ in particular, so $F(z)$ has a Laurent expansion in powers of $e^{2\pi i z}$. The condition of moderate growth rules out negative powers. We may and shall exclude even the 0-th power – in other words, we shall assume $F$ is cuspidal. Eisenstein series, which span a linear complement of the space of cuspidal modular forms, require a separate but simpler analysis. In the cuspidal case, then,

$$
\text{(3.10)} \quad F(z) = \sum_{n=1}^{\infty} a_n\, n^{(k-1)/2}\, e^{2\pi i n z}.
$$

Traditionally one denotes $e^{2\pi i z}$ by the symbol $q$ and calls (3.9) the $q$-expansion, but we shall not follow that practice.

The Dirichlet series with coefficients $a_n$ is called the standard $L$-function of the modular form $F$,

$$
\text{(3.11)} \quad L(s, F) = \sum_{n=1}^{\infty} a_n\, n^{-s}.
$$

It satisfies a functional equation, which reflects the transformation law for $F$ with respect to $z \mapsto -1/z$. We should remark that the factor $|n|^{(k-1)/2}$ in (3.10) has the effect of making the functional equation relate $L(s, F)$ to $L(1-s, F)$ rather than to $L(k-s, F)$, as is the case in the classical literature. The $L$-function (3.11) can be "twisted", either by an additive character $\chi$ of finite order,

$$
\text{(3.12)} \quad \chi(n) = e^{2\pi i n a/c}, \text{ with } a, c \in \mathbb{Z},\ c \neq 0,\ (a,c) = 1,
$$

or by a Dirichlet character modulo $q$,

$$
\text{(3.13)} \quad \chi(n) = \begin{cases} \chi_q([n]) & \text{if } (n,q) = 1 \\ 0 & \text{if } (n,q) \neq 1, \end{cases}
$$
$$
\text{for some character } \chi_q : (\mathbb{Z}/q\mathbb{Z})^* \to \mathbb{C}^*,
$$

with $[n]$ denoting the image of $n$ in $(\mathbb{Z}/q\mathbb{Z})^*$. The Dirichlet series

$$
\text{(3.14)} \quad L_\chi(s, F) = \sum_{n=1}^{\infty} \chi(n)\, a_n\, n^{-s}
$$

is the twisted $L$-function of $F$, with "additive twist" in the case of (3.12), respectively "multiplicative twist" in the case of (3.13). Multiplicative twists have arithmetic significance. For example, multiplicatively twisted $L$-functions of Hecke eigenforms have Euler products and are considered special cases of Langlands $L$-functions. Additive twists, on the other hand, play an important role in the analytic study of $L$-functions.

A Dirichlet character $\chi$ modulo $q$ is called primitive if, in the notation of (3.13), $\chi_q$ does not arise as the pullback of a character of $(\mathbb{Z}/q'\mathbb{Z})^*$ for some divisor $q'$ of $q$. In this situation, the multiplicatively twisted $L$-function $L_\chi(s,F)$ also satisfies a functional equation. Everything that has been said so far can be generalized in some fashion to automorphic forms on $GL(n)$ and other higher rank groups. Not so the assertion that the additively twisted $L$-function $L_\chi(s,F)$ satisfies a functional equation, which is true for $GL(2)$ but fails already for $GL(3)$.

Maass forms are eigenfunctions of the Laplace operator on the upper half plane, of moderate growth, invariant under the natural action of $SL(2,\mathbb{Z})$ or more generally, invariant under a subgroup $\Gamma \subset SL(2,\mathbb{Z})$ of finite index. Just as in the case of modular forms, one can attach $L$ functions to Maass forms. These $L$-functions, both with multiplicative and additive twists, again satisfy functional equations.

In 1930, Wilton proved a Voronoi summation formula for the modular form $\Delta$, the essentially unique cusp form of weight 12. Wilton's argument applies to any cuspidal modular form, even with additive twist. It can also be adapted to Maass forms. The reason, of course, is the functional equation for the additively twisted $L$-functions, which is used in the same way as the functional equations for $\zeta(s)$, $\zeta^2(s)$, and $\zeta_{\mathbb{Q}(\sqrt{-1})}(s)$. This type of argument does not extend to groups beyond $GL(2)$, however, since then additively twisted $L$-functions no longer satisfy functional equations. Recently we were able to prove a Voronoi formula for $GL(3)$ with representation-theoretic arguments [19]. It seems clear that our argument works for $GL(n)$.

In the next section we state the Voronoi summation formula for modular forms and Maass forms and then sketch a proof. Our argument is essentially different from the usual arguments, in that it works directly with automorphic forms; in effect, it is an adaptation to $GL(2)$ of our argument for $GL(3)$. We go on to discuss our $GL(3)$-formula in section 5.

## 4 Voronoi Summation for $GL(2)$

Cuspidal Maass forms and cuspidal modular forms for $\Gamma = SL(2,\mathbb{Z})$ correspond to embeddings of, respectively, principal series and discrete series representations into $L^2(\Gamma \backslash SL(2,\mathbb{R}))$. Let us begin with some quite general remarks about discrete summands in $L^2(\Gamma \backslash G)$.

Initially $G$ will denote an arbitrary unimodular Lie group of type I, $\Gamma \subset G$ a discrete subgroup, and $\omega : Z_G \to \mathbb{C}^*$ a unitary character of $Z_G =$ center of $G$.

Then $G$ acts unitarily on

$$
(4.1) \qquad \begin{aligned} L^2_\omega(\Gamma\backslash G) \ = \ \{\, f \in L^2_{\text{loc}}(\Gamma\backslash G) \mid f(gz) = \omega(z)\, f(g) \\ \text{for } z \in Z_G,\ \text{and}\ \ \textstyle\int_{\Gamma\backslash G/Z_G} |f|^2\, dg < \infty\,\}. \end{aligned}
$$

If $(\pi, V)$ is an irreducible unitary representation of $G$ and

$$
(4.2) \qquad j : V \ \hookrightarrow \ L^2_\omega(\Gamma\backslash G)
$$

a $G$-invariant embedding, the subspace $V^\infty \subset V$ of $C^\infty$ vectors gets mapped to $\Gamma$-invariant smooth functions on $G$, which can be evaluated at the identity. Thus

$$
(4.3) \qquad \tau_j : V^\infty \ \longrightarrow\ \mathbb{C}\,, \qquad \tau_j(v)\ =\ j(v)(e)\,,
$$

is a well defined, $\Gamma$-invariant linear map. Since $j(v)(e)$ can be bounded in terms of finitely many $L^2$ derivatives, one can show that $\tau_j$ is continuous with respect to the natural topology on $V^\infty$. In short, $\tau_j$ can be regarded as a $\Gamma$-invariant distribution vector for the unitary representation $(\pi', V')$ dual to $(\pi, V)$,

$$
(4.4) \qquad \tau_j\ \in (V'^{-\infty})^\Gamma\,.
$$

We shall call $\tau_j$ the automorphic distribution corresponding to the embedding $j$. It determines $j$ completely since $j(v)(g) = j(\pi(g)v)(e) = \tau_j(\pi(g)v)$ for any $v$ in the dense subspace $V^\infty \subset V$. Our discussion will focus on $\tau_j$. To simplify the notation, we switch the roles of $(\pi, V)$ and $(\pi', V')$; from now on,

$$
(4.5) \qquad j : V' \ \hookrightarrow\ L^2_{\omega^{-1}}(\Gamma\backslash G)\,,\qquad \tau_j\ \in (V^{-\infty})^\Gamma\,.
$$

This is legitimate since integration over $\Gamma\backslash G/Z_G$ sets up a nondegenerate pairing between $L^2_\omega(\Gamma\backslash G)$ and $L^2_{\omega^{-1}}(\Gamma\backslash G)$.

We now specialize our hypotheses to the case of $G = SL(2,\mathbb{R})$ and $\Gamma = SL(2,\mathbb{Z})$. Since $\Gamma$ contains $Z_G$, only discrete series representations with even minimal weight and even principal series representations can occur as summands in (4.4). Discrete series representations have realizations as subrepresentations of non-unitary principal series representations. According to a result of Casselman and Wallach [1] – whose specialization to $G = SL(2,\mathbb{R})$ can be verified directly – the assignment $V \mapsto V^{-\infty}$ is an exact functor. Thus, even when $\pi$ is a representation of the discrete series, the automorphic distribution $\tau_j$ may be regarded as a distribution vector for a (possibly reducible) even principal series representation. As model for even principal series representations we take

$$
(4.6) \qquad V^\infty_\nu\ =\ \{\, f \in C^\infty(\mathbb{R})\ \mid\ |x|^{2\nu-1} f(1/x) \in C^\infty(\mathbb{R})\,\}\,,
$$

with $\nu \in \mathbb{C}$ as representation parameter. The group $G$ acts via

$$
(4.7) \qquad (\pi_\nu(g)f)(x)\ =\ |cx+d|^{2\nu-1} f\!\left(\tfrac{ax+b}{cx+d}\right) \quad \text{if}\quad g^{-1} = \begin{pmatrix} a & b \\ c & d \end{pmatrix}.
$$

The resulting representation is irreducible whenever $\nu \notin \mathbb{Z} + 1/2$; in that case $\pi_\nu \simeq \pi_{-\nu}$. If $\nu = 1/2 - k \in \mathbb{Z} + 1/2$ is negative, $V_\nu^\infty$ contains two discrete series representations as subrepresentations and has a $(2k-1)$-dimensional quotient, whereas $V_{-\nu}^\infty$ has a $(2k-1)$-dimensional subrepresentation and two discrete series quotients.

In contrast to (4.6), the restriction map from the space of distribution vectors $V_\nu^{-\infty}$ to $C^{-\infty}(\mathbb{R})$ is not injective: a distribution on a manifold is not determined by its restriction to a dense open subset. Still, in terms of the restriction map

$$(4.8) \qquad V_\nu^{-\infty} \longrightarrow C^{-\infty}(\mathbb{R}),$$

the action of $G$ is given by the formula (4.7) at all points $x \in \mathbb{R}$ such that $cx + d \neq 0$.

We now suppose that $\tau_j \in (V_\nu^{-\infty})^\Gamma$ is a cuspidal automorphic distribution – i.e., an automorphic distribution corresponding to a cuspidal embedding (4.5). The automorphy with respect to $\Gamma$ implies in particular that $\tau_j$ is periodic of period 1. Consequently its restriction to $\mathbb{R}$ has a Fourier expansion

$$(4.9) \qquad \tau_j(x) = \sum_{n \neq 0} a_n \, |n|^{-\nu} \, e^{2\pi i n x},$$

whose constant term vanishes because of the assumption of cuspidality. In the case of an embedding of a discrete series representation, $a_n = 0$ for either all positive or all negative indices $n$. In the case of a principal series representation, the parameter $\nu$ is determined only up to sign. The two possible choices determine different distributions in $C^{-\infty}(\mathbb{R})$, which are related by the standard intertwining operator. One can show that the coefficients $a_n$ corresponding to the two choices agree up to a multiplicative constant – that is accomplished by the factor $|n|^{-\nu}$ in (4.9). If $\tau_j$ satisfies the conditions of automorphy under $\Gamma$ and cuspidality, then so does $x \mapsto \tau_j(-x)$. For discrete series representations, this operation switches the roles of the two subrepresentations in the ambient $V_\nu^{-\infty}$. In the Maass case, on the other hand, this operation allows us to consider even and odd automorphic distributions separately. In other words, in the Maass case we may and shall assume

$$(4.10) \qquad \tau_j(-x) = (-1)^\eta \tau_j(x) \qquad (\eta \in \mathbb{Z}/2\mathbb{Z}).$$

When this condition is imposed, the $a_n$ for $n > 0$ determine also the $a_n$ with negative index, and

$$(4.11) \qquad L(s, \tau_j) = \sum_{n=1}^\infty a_n \, n^{-s}$$

is the standard $L$-function of the cuspidal modular form or Maass form corresponding to the embedding (4.5); one sees this most easily by identifying the action of the Hecke operators on the Fourier series (4.9).

For the statement of the Voronoi formula for $GL(2)$, we fix an automorphic distribution $\tau_j$ as above, as well as integers $a, \bar{a}, c$ such that $c \neq 0$, $(a, c) = 1$, and $\bar{a} a \equiv 1$ modulo $c$.

**4.12 Theorem.** [2, 11, 15] *If the Schwartz function $f \in \mathcal{S}(\mathbb{R})$ vanishes to infinite order at the origin,*

$$\sum_{n \neq 0} a_n\, e^{-2\pi i n a/c}\, f(n) \;=\; |c| \sum_{n \neq 0} \frac{a_n}{|n|}\, e^{2\pi i n \bar{a}/c}\, F(n/c^2),$$

*with $F$ described by the conditionally convergent integral*

$$F(t) = \iint_{\mathbb{R}^2} f(x_1 x_2/t)\, |x_1|^\nu\, |x_2|^{-\nu}\, e^{-2\pi i (x_1+x_2)}\, dx_1\, dx_2\,.$$

Care must be taken in interpreting the integral expression for $F$. When the integration is performed in the indicated order – first with respect to $x_1$, then with respect to $x_2$ – the inner integral converges absolutely, but the outer integral only conditionally, and only when $\operatorname{Re} \nu > 0$. For other values of $\nu$, the integral retains meaning by analytic continuation. Alternatively $F$ can be characterized by a formula for its Mellin transform in terms of the Mellin transform of $f$; this is the form in which the theorem is usually stated. The proof shows that $F$ and all of its derivatives decay rapidly as $|x| \to \infty$. At the origin $F$ has certain mild singularities; see [19] for details.

The Voronoi summation formula for $GL(2)$ has become an essential tool in analytic number theory, for problems such as sub-convexity of automorphic $L$-functions or the analysis of the Petersson and Kuznetsov trace formulas. In all the applications the presence of the additive twist plays a critical role. A formal application of the theorem with $f(x) = (\operatorname{sgn} x)^\eta |x|^{-s}$ – which is not covered by the hypotheses, of course – gives the functional equation for the $L$-function with additive twists. In fact, the proof that we shall sketch can be adapted to deal even with this "illegitimate" choice of test function. The functional equation for multiplicative twists can be deduced by taking appropriate linear combinations of the additively twisted functional equation. Conversely, the functional equation with multiplicative twist implies the equation with additive twist, and that in turn can be used to establish the summation formula, along the lines of the arguments in section 2 [2, 11, 15].

For $GL(3)$, the passage from multiplicative to additive twists is far more subtle; in fact, additively twisted $L$-functions for $GL(3)$ do not satisfy functional equations, as was pointed out already. This has been the obstacle to a Voronoi formula for $GL(3)$ in the past. Our proof of the $GL(3)$ formula in [19] overcomes the obstacle by working directly with the automorphic distribution. We shall now sketch a proof of theorem 4.12 by adapting our argument in [19] to the case of $SL(2, \mathbb{R})$.

Because of the hypotheses of the theorem, there exits $b \in \mathbb{Z}$ such that $a\bar{a} = 1 + bc$. The $2 \times 2$ matrix $\gamma$ with entries $\bar{a}$, $b$, $c$, $a$ lies in $\Gamma$. We now use the formula (4.7) to make the automorphy of $\tau_j$ under $\gamma$ explicit, then substitute $x - a/c$ for $x$, which results in the equation

(4.13) $\qquad \tau_j\left(x - \tfrac{a}{c}\right) \;=\; |cx|^{2\nu-1}\, \tau_j\left(\tfrac{\bar{a}}{c} - \tfrac{1}{c^2 x}\right).$

We integrate this equality of tempered distributions against a test function $g$ in the Schwartz space $\mathcal{S}(\mathbb{R})$. On one side, using the Fourier expansion (4.9), we get

$$
\begin{aligned}
\int_{\mathbb{R}} \tau_j \left(x - \tfrac{a}{c}\right) g(x)\, dx \ &= \\
(4.14) \qquad &= \int_{\mathbb{R}} \sum_{n \neq 0} a_n\, |n|^{-\nu}\, e^{2\pi i n(x - a/c)}\, g(x)\, dx \\
&= \sum_{n \neq 0} a_n\, |n|^{-\nu}\, e^{-2\pi i n a/c}\, \widehat{g}(-n)\,.
\end{aligned}
$$

The other side of (4.13) appears to be ill defined at $x = 0$. However, the fact that $\tau_j$ is periodic without constant term implies that this distribution can be expressed as the $N$-th derivative of a bounded continuous function for every sufficiently large integer $N$. This fact, coupled with an integration-by-parts argument, can be used to show that the distribution $\tau_j(1/x)$ has a canonical extension across the origin. With this interpretation, (4.13) remains valid even at $x = 0$; roughly speaking, the cuspidality hypothesis rules out the delta function or derivatives of the delta function from contributing to the expression for of $\tau_j$ at infinity. In short, the formal computation

$$
\begin{aligned}
\int_{\mathbb{R}} |cx|^{2\nu - 1}\, \tau_j\!\left(\tfrac{\bar{a}}{c} - \tfrac{1}{c^2 x}\right) g(x)\, dx \ &= \\
(4.15) \qquad &= |c|^{2\nu - 1} \int_{\mathbb{R}} |x|^{2\nu - 1} \sum_{n \neq 0} a_n\, |n|^{-\nu}\, e^{2\pi i n c^{-2}(\bar{a} c - 1/x)}\, g(x)\, dx \\
&= |c|^{2\nu - 1} \sum_{n \neq 0} a_n\, |n|^{-\nu} \int_{\mathbb{R}} |x|^{2\nu - 1}\, e^{2\pi i n c^{-2}(\bar{a} c - 1/x)}\, g(x)\, dx
\end{aligned}
$$

can be made rigorous [20].

Let $f \in \mathcal{S}(\mathbb{R})$ be a Schwartz function which vanishes to infinite order at $x = 0$, as in the statement of the theorem. Then $g(x) = \int_{\mathbb{R}} f(y)|y|^\nu e(-xy)\, dy$ is also a Schwartz function, and $f(x) = |x|^{-\nu}\, \widehat{g}(-x)$. With this choice of $g$, (4.13–4.15) imply

$$
\begin{aligned}
\sum_{n \neq 0} a_n\, e^{-2\pi i n a/c}\, f(n) \ &= \\
&= \sum_{n \neq 0} a_n\, \frac{|c|^{2\nu - 1}}{|n|^\nu} \iint_{\mathbb{R}^2} f(y)\, |y|^\nu\, |x|^{2\nu - 1}\, e^{2\pi i(n\bar{a}/c - n/c^2 x - xy)}\, dy\, dx \\
&= \sum_{n \neq 0} a_n\, \frac{|c|^{2\nu - 1}}{|n|^\nu} \iint_{\mathbb{R}^2} f(y)\, |y|^\nu\, |x|^{-2\nu - 1}\, e^{2\pi i(n\bar{a}/c - n/c^2 x - y/x)}\, dy\, dx \\
&= \sum_{n \neq 0} a_n\, \frac{|c|^{2\nu - 1}}{|n|^\nu} \iint_{\mathbb{R}^2} f(xy)\, |y|^\nu\, |x|^{-\nu}\, e^{2\pi i(n\bar{a}/c - n/c^2 x - y)}\, dy\, dx
\end{aligned}
$$

$$= \sum_{n \neq 0} a_n \frac{|c|}{|n|} \iint_{\mathbb{R}^2} f(c^2 xy/n) \, |y|^\nu \, |x|^{-\nu} \, e^{2\pi i(n\bar{a}/c - x - y)} \, dy \, dx$$

$$= |c| \sum_{n \neq 0} \frac{a_n}{|n|} e^{2\pi i n\bar{a}/c} \iint_{\mathbb{R}^2} f(c^2 xy/n) \, |y|^\nu \, |x|^{-\nu} e^{-2\pi i(x+y)} \, dy \, dx \, .$$

Setting $t = nc^{-2}$ and changing the symbols for the variables of integration we get the statement of the theorem.

## 5 The $GL(3)$ summation formula

The discussion in the beginning of the previous section applies also to the setting of $G = GL(3, \mathbb{R})$, $\Gamma = GL(3, \mathbb{Z})$. Now $G$ has one dimensional center, so the central character $\omega$ in (4.1–4.5) is necessary – without it, the space of $\Gamma$-invariant $L^2$ functions would vanish.

As in the case of $SL(2, \mathbb{R})$, every cuspidal automorphic distribution can be regarded as a distribution vector for some principal series representation; this assertion depends on the Casselman embedding theorem and the Casselman-Wallach exactness theorem for spaces of distribution vectors. Not-necessarily-unitary principal series representations of $G$ with unitary central character correspond to parameters $(\lambda, \delta)$,

(5.1) $$\begin{aligned} \lambda &= (\lambda_1, \lambda_2, \lambda_3) \in \mathbb{C}^3, \quad \delta = (\delta_1, \delta_2, \delta_3) \in (\mathbb{Z}/2\mathbb{Z})^3, \\ &\sum_j \operatorname{Re} \lambda_j = 0 \, . \end{aligned}$$

The $\delta_j$ specify the inducing character on the component group of the diagonal Cartan subgroup, and the $\lambda_j$ are normalized so that the Weyl group action fixes $\lambda = 0$ and $\rho = (1, 0, -1)$ is the half-sum of the positive roots. The condition $\sum_j \operatorname{Re} \lambda_j = 0$ reflects the requirement that $Z_G$ must act unitarily. As model for $V_{\lambda,\delta}^\infty$ we take the space of smooth functions on

(5.2) $$N_+ = \left\{ \begin{pmatrix} 1 & x & z \\ 0 & 1 & y \\ 0 & 0 & 1 \end{pmatrix} \, \middle| \, x, y, z \in \mathbb{R} \right\} \subset GL(3, \mathbb{R})$$

with an appropriate regularity condition at "infinity" and an action $\pi_{\lambda,\delta}$ of $G$ involving a factor of automorphy. In analogy to (4.8), there is a natural restriction map

(5.3) $$V_{\lambda,\delta}^{-\infty} \longrightarrow C^{-\infty}(N_+) \, .$$

The subgroup $N_+ \subset G$ acts by translation, without factor of automorphy, via this restriction map.

We now consider a particular cuspidal automorphic distribution $\tau_j \in V_{\lambda,\delta}^{-\infty}$. By restriction to $N_+$, and in terms of the coordinates $x, y, z$ in (5.2), we may think of $\tau_j$ as a distribution in these variables. Naively speaking, the restriction to $N_+$ does not determine $\tau_j$, but the $\Gamma$-automorphy of $\tau_j$ ensures that no information is lost. The invariance under $N_+ \cap \Gamma$ allows us to regard $\tau_j(x,y,z)$ as a distribution on the compact manifold $(N_+ \cap \Gamma)\backslash N_+$. As such, it has a Fourier decomposition with respect to the right action of $N_+$:

$$
(5.4) \quad \begin{aligned} \tau_j(x,y,z) &= \\ &= \sum_{m,n\neq 0} \frac{|n|^{\lambda_3}(\operatorname{sgn} n)^{\delta_3}}{|m|^{\lambda_1}(\operatorname{sgn} m)^{\delta_1}}\, a_{m,n}\, e^{2\pi i(mx+ny)} \;+\; \ldots\,; \end{aligned}
$$

here we have written out the abelian part of the decomposition explicitly, but the remaining components – those on which the center of $N_+$ acts according to a non-trivial character – are summarily denoted by $\ldots$. The cuspidality of $\tau_j$ accounts for the absence of abelian terms indexed by $(0,n)$ or $(m,0)$. The powers of $|m|$, $|n|$ and the signs play the same role as the factor $|n|^{-\nu}$ in (4.9) – they make the coefficients independent of the choice of Casselman embedding, up to universal multiplicative factors. With this normalization of the coefficients,

$$
(5.5) \quad L(s,\tau_j) \;=\; \sum_{n=1}^{\infty} a_{1,n}\, n^{-s}\,, \qquad L(s,\tilde\tau_j) \;=\; \sum_{n=1}^{\infty} a_{n,1}\, n^{-s}
$$

are the standard $L$-functions of the automorphic representation corresponding to $\tau_j$ and the dual automorphic representation; the latter corresponds to the dual automorphic distribution $\tilde\tau_j$ whose abelian Fourier coefficients agree with those of $\tau_j$ except for a transposition of the indices.

The statement of our Voronoi formula for $GL(3)$ is entirely analogous to that for $GL(2)$. It involves the Kloosterman sum

$$
(5.6) \quad S(m,n;c) \;=\; \sum_{x\in(\mathbb{Z}/c\mathbb{Z})^*} e^{2\pi i(mx+n\bar x)/c} \quad (m,n,c\in\mathbb{Z},\; c\neq 0)\,;
$$

in the summation, $\bar x$ denotes an integer which represents the multiplicative inverse of $x$ modulo $c$.

**5.7 Theorem.** *Let $f \in \mathcal{S}(\mathbb{R})$ be a Schwartz function which vanishes to infinite order at the origin, or more generally, a function on $\mathbb{R} - \{0\}$ such that $(\operatorname{sgn} x)^{\delta_3}|x|^{-\lambda_3} f(x) \in \mathcal{S}(\mathbb{R})$. Then for $(a,c)=1$, $c\neq 0$, $\bar a a \equiv 1 \pmod{c}$ and $q > 0$,*

$$
\begin{aligned}
\sum_{n\neq 0} a_{q,n}\, e^{-2\pi i n a/c} f(n) &= \\
&= \sum_{d|cq} \left|\frac{c}{d}\right|^{1-\lambda_1-\lambda_2-\lambda_3} \sum_{n\neq 0} \frac{a_{n,d}}{|n|}\, S(q\bar a, n; qc/d)\, F\!\left(\frac{nd^2}{c^3 q}\right),
\end{aligned}
$$

*where*

$$
F(t) \;=\; \int_{\mathbb{R}^3} f\!\left(\frac{x_1 x_2 x_3}{t}\right) \prod_{j=1}^{3} \left(e^{-2\pi i x_j} |x_j|^{-\lambda_j}(\operatorname{sgn} x_j)^{\delta_j}\right) dx_3\, dx_2\, dx_1.
$$

The integral defining $F$ converges when performed as repeated integral in the indicated order – i.e., with $x_3$ first, then $x_2$, then $x_1$ – and provided $\operatorname{Re}\lambda_1 > \operatorname{Re}\lambda_2 > \operatorname{Re}\lambda_3$; it has meaning for arbitrary values of $\lambda_1$, $\lambda_2$, $\lambda_3$ by analytic continuation. As in the case of $GL(2)$, the function $F$ and all its derivatives decay rapidly as $|x| \to \infty$, but $F$ has a mild singularity at the origin. One can also characterize $F$ in terms of its Mellin transform [19].

While the statement of the theorem only involves the abelian Fourier coefficients, the proof uses the non-abelian coefficients in an essential way. The invariance of $\tau_j$ under the two simple Weyl reflections – which have representatives in $\Gamma$ – translates into relations between the abelian and non-abelian coefficients. Additive twists show up transparently in these relations. A careful analysis then leads to the summation formula. Variants of these arguments can be used to prove the functional equations for the standard $L$-functions, with multiplicative twists, and also to prove the "converse theorem" of Jacquet, Piatetski-Shapiro and Shalika [14]: an automorphic representation can be reconstructed from the two standard $L$-functions, provided they have Euler products of the appropriate type, their multiplicative twists satisfy functional equations, and they grow slowly in the imaginary directions [19].

Cuspidal $GL(3,\mathbb{Z})$-automorphic representations exist abundantly, as a consequence of variants of the Selberg trace formula [18]. Explicit construction of such representations is a very different matter – for that, one relies mainly on the Gelbart-Jacquet symmetric square lift of cuspidal $SL_2(\mathbb{Z})$-automorphic representations [4]. Not surprisingly, then, our $GL(3)$ formula has been applied first to a problem about modular forms on $GL(2)$. Indeed, when Peter Sarnak saw our analysis of automorphic distributions on $GL(3,\mathbb{R})$, he suggested we should be able to establish a Voronoi summation formula for $GL(3)$, a formula he considered potentially useful for the solution of a certain problem about Maass forms. Let us describe this problem.

Let $X$ be a compact Riemann surface, $\{\phi_j\}_{j=0}^\infty$ an orthonormal basis of Laplace eigenfunctions, and $\mu_j$ the eigenvalue corresponding to the $j$-th eigenfunction. Sogge [26] bounds the $L^p$-norms of the $\phi_j$ in terms of the $\mu_j$. His bound comes from differential-geometric estimates and is sharp when applied to a general compact Riemann surface. A stronger bound is expected to hold for the discrete spectrum of a complete hyperbolic surface of finite volume, even when the surface fails to be compact. In the special case of the surface $X = SL_2(\mathbb{Z})\backslash\mathcal{H}$, $\mathcal{H} = \{\operatorname{Re} z > 0\}$, Berry's random wave model and experimental data of Hejhal-Rackner [9] suggest that the distribution of the values of eigenfunctions tend towards a Gaussian distribution. In particular this would imply the estimate

$$(5.8) \qquad \|\phi_j\|_p \;=\; O(\mu_j^\epsilon) \qquad (\, 2 < p < \infty, \;\; \epsilon > 0 \,)\,,$$

for any orthonormal basis of Hecke-Laplace eigenfunctions of the discrete spectrum of the modular surface $SL_2(\mathbb{Z})\backslash\mathcal{H}$. The survey [23] explains these conjectures as eigenfunction analogues of the Lindelöf conjecture; in fact, a variant of (5.8) would imply the classical Lindelöf conjecture.

Sarnak and Watson [24] have recently announced (5.8) for the case $p = 4$ – at present under the assumption of the Ramanujan conjectures for Maass forms on $GL(2)$. In contrast, the bound of [26] for $p = 4$, were it to apply in this noncompact situation, only predicts the bound $O(\mu_j^{1/16})$. Sarnak-Watson's argument uses our Voronoi summation formula for $GL(3)$, among other ingredients.

To show how our formula can be applied to modular forms and Maass forms, we shall describe the abelian Fourier coefficients of the symmetric square lifting terms of the coefficients of the original modular or Maass form. Thus we consider a cuspidal $SL(2, \mathbb{Z})$-automorphic distribution $\tau$, with representation parameter $\nu$, which we express as in (4.9). Recall that $a_n = 0$ for $n < 0$ and $\nu = 1/2 - k$ if $\tau$ corresponds to a holomorphic cusp form of weight $2k$. In the Maass case, we impose the parity condition (4.10), which is equivalent to

$$(5.9) \qquad a_{-n} = (-1)^\eta a_n \qquad (\eta \in \mathbb{Z}/2\mathbb{Z}).$$

Still in the Maass case, the parameter $\nu$ is purely imaginary, but only determined up to a sign factor. In both cases, the $a_n$ with $n > 0$ completely determine $\tau$.

The symmetric square lift is defined for Hecke eigenforms. The Hecke property implies in particular that $a_1 \neq 0$, making it possible to rescale $\tau$ so that $a_1 = 1$. In that case, the Hecke property is equivalent to the Euler product factorization

$$(5.10) \qquad L(s, \tau) = \sum_{n=1}^\infty a_n n^{-s} = \prod_p (1 - \alpha_p p^{-s})^{-1} (1 - \alpha_p^{-1} p^{-s})^{-1}$$

$$\text{with} \quad \alpha_p + \alpha_p^{-1} = a_p.$$

The Gelbart-Jacquet symmetric square lift of $\tau$ is a cuspidal $GL(3, \mathbb{Z})$-automorphic distribution $\mathrm{Sym}^2 \tau$, whose abelian Fourier coefficients we denote by $a_{m,n}$. Recall the definition of the Möbius function $\mu : \mathbb{N} \to \{\pm 1\}$,

$$(5.11) \qquad \mu(p_1^{e_1} p_2^{e_2} \cdots p_r^{e_r}) = \begin{cases} (-1)^r & \text{if all } e_j = 1 \text{ and} \\ & \quad p_i \neq p_j \text{ for } i \neq j, \\ 0 & \text{if at least one } e_j > 1. \end{cases}$$

**5.12 Proposition.** *The lifted automorphic distribution $\mathrm{Sym}^2 \tau$ has representation parameters*

$$\lambda = (2\nu, -2\nu, 0), \quad \delta = \begin{cases} (0, 0, 0) & \text{in the Maass case} \\ (1, 0, 1) & \text{in the holomorphic case} \end{cases}$$

*and abelian Fourier coefficients*

$$a_{m,n} = \sum_{d \mid (m,n)} \mu(d) \sum \left\{ a_r a_s \mid r = \frac{m^2}{k^4 d^2}, \ s = \frac{n^2}{\ell^4 d^2}, \ k, \ell \in \mathbb{N} \right\}.$$

With these choices of $\lambda$, $\delta$ and $a_{m,n}$, theorem 5.7 provides information about the $a_n$. In the Maass case, the representation $\pi_{\lambda,\delta}$ of $GL(3,\mathbb{R})$ which corresponds to the parameters $(\lambda, \delta)$ is an irreducible principle series representation, and in the holomorphic case, $\pi_{\lambda,\delta}$ is unitarily and irreducibly induced from a discrete series representation of $GL(2,\mathbb{R})$.

*Proof.* The completed $L$-function of the automorphic representation determined by $\tau$ is the product

$$(5.13) \qquad \Lambda(s,\tau) \;=\; L_\infty(s,\tau)\, L_{\text{finite}}(s,\tau)\,,$$

with $L_{\text{finite}}(s,\tau) = L(s,\tau)$, as in (5.10), corresponding to the finite places and

$$(5.14) \qquad \begin{aligned} L_\infty(s,\tau) &= \\ &= \begin{cases} \Gamma_\mathbb{R}(s+\nu+\eta)\,\Gamma_\mathbb{R}(s-\nu+\eta) & \text{in the Maass case} \\ \Gamma_\mathbb{C}(s+k-\tfrac{1}{2}) & \text{in the holomorphic case} \end{cases} \end{aligned}$$

corresponding to the archimedean place. Here we are following the notational convention

$$(5.15) \qquad \Gamma_\mathbb{R}(s) \;=\; \pi^{-s/2}\,\Gamma(s/2)\,, \qquad \Gamma_\mathbb{C}(s) \;=\; 2\,(2\pi)^{-s}\,\Gamma(s)\,,$$

which is commonly used in the context of $L$-functions. The completed $L$-function (5.13) satisfies the functional equation

$$(5.16) \qquad \Lambda(s,\tau) \;=\; \omega\,\Lambda(1-s,\tau)\,, \qquad \text{with } \omega \in \mathbb{C}\,,\ |\omega|=1\,.$$

Gelbart-Jacquet [4] characterize the symmetric square lifting $\operatorname{Sym}^2 \tau$ by its completed $L$-function

$$(5.17) \qquad \Lambda(s, \operatorname{Sym}^2 \tau) \;=\; L_\infty(s, \operatorname{Sym}^2 \tau)\, L_{\text{finite}}(s, \operatorname{Sym}^2 \tau)\,,$$

which satisfies the same type of functional equation as $\Lambda(s,\tau)$. The archimedean part is

$$(5.18) \qquad \begin{aligned} L_\infty(s, \operatorname{Sym}^2 \tau) &= \\ &= \begin{cases} \Gamma_\mathbb{R}(s+2\nu)\Gamma_\mathbb{R}(s)\Gamma_\mathbb{R}(s-2\nu) & \text{in the Maass case} \\ \Gamma_\mathbb{C}(s+2k-1)\Gamma_\mathbb{R}(s+1) & \text{in the holomorphic case} \end{cases} \end{aligned}$$

and the finite part

$$(5.19) \qquad \begin{aligned} L_{\text{finite}}(s, \operatorname{Sym}^2 \tau) &= \sum_{n=1}^\infty A_n\, n^{-s} \;=\; \\ &= \prod_p (1-\alpha_p^2 p^{-s})^{-1}\,(1-p^{-s})^{-1}\,(1-\alpha_p^{-2} p^{-s})^{-1}. \end{aligned}$$

Expanding the Euler product one finds

$$(5.20) \qquad L_{\text{finite}}(s, \operatorname{Sym}^2 \tau) \;=\; \zeta(2s) \sum_{n=1}^\infty a_{n^2}\, n^{-s}\,,$$

which in turn implies

$$(5.21) \qquad A_n \;=\; \sum\nolimits_{d^2 \mid n} a_{(n/d^2)^2}\,.$$

In [19, Theorem 6.15] we relate the coefficients $a_{n,m}$ and parameters $\lambda, \delta$ of a cuspidal $GL(3,\mathbb{Z})$-automorphic representation of $GL(3,\mathbb{R})$ to its standard $L$-function. In the current setting this gives

$$(5.22) \qquad a_{1,n} \;=\; a_{n,1} \;=\; A_{|n|} \qquad (n \neq 0)\,.$$

The abelian Fourier coefficients $a_{m,n}$ of any Hecke eigendistribution for $GL(3,\mathbb{Z})$ satisfy the product relation

$$(5.23) \qquad a_{m,n} \;=\; \sum\nolimits_{d \mid (m,n)} \mu(d)\, a_{1,n/d}\, a_{m/d,1}\,.$$

In the situation at hand this leads to the formula for the $a_{m,n}$ asserted by the proposition. The values of the parameters $\lambda, \delta$ in the statement of the proposition match (5.18) to the calculation of the $\Gamma$-factors in [19, Theorem 6.15].

## 6  Error term in the circle problem

In this section we try to give a glimpse of how the Voronoi summation formula is used in practice. Voronoi proved the bound

$$(6.1) \qquad \Delta(X) \;=\; O(X^{1/3})$$

for the error term (1.5). We shall show how to obtain the slightly weaker bound

$$(6.2) \qquad \Delta(X) \;=\; O(X^{1/3+\epsilon})$$

from fairly simple estimates. Voronoi's formula (1.7) for the sequence $r_2(n)$ can be established from the Poisson summation formula in two variables, applied to radially symmetric test functions $f$. In the interest of brevity we shall deduce (1.7) directly from Poisson summation in two variables, instead of taking a detour via (1.7) – which would not fundamentally change the argument.

Let $\chi_D$ denote the characteristic function of $D$, the unit disc in $\mathbb{R}^2$. In the following, $X > 1$ will be a large parameter and $\delta = \delta(X) = X^{-1/6}$. We choose a radially symmetric, smooth function $\Phi \geq 0$, supported on the unit disc $D$, of total integral 1. The convolution product

$$(6.3) \qquad F_X(v) \;=\; \delta^{-2}\, \chi_D(X^{-1/2}v) * \Phi(\delta^{-1}v)$$

is also radially symmetric and smooth, $0 \leq F_X \leq 1$, $F_X(v) = 1$ for $\|v\| \leq \sqrt{X} - \delta$ and $F_X(v) = 0$ for $\|v\| \geq \sqrt{X} + \delta$. Its Fourier transform

$$(6.4) \qquad \widehat{F_X}(v) \;=\; X\, \widehat{\chi_D}(X^{1/2}v)\, \widehat{\Phi}(\delta v)$$

decays rapidly and is again radially symmetric. Because of the radial symmetry, there exist functions $f_X, \psi : \mathbb{R}_{\geq 0} \to \mathbb{R}$ such that

(6.5) $$F_X(v) = f_X(\|v\|), \qquad \widehat{\Phi}(v) = \psi(\|v\|).$$

The Fourier transform of $\widehat{\chi_D}$ is expressible in terms of a Bessel function:

(6.6) $$\widehat{\chi_D}(v) = \|v\|^{-1} J_1(2\pi\|v\|)$$

[6, p. 962]. Since $\chi_D$ and $\Phi$ have total integral $\pi$ and 1, respectively,

(6.7) $$\widehat{\chi_D}(0) = \pi, \qquad \psi(0) = 1;$$

the former also follows from the known properties of $J_1$, of course. Poisson summation in $\mathbb{R}^2$ identifies the sum $\sum_{v \in \mathbb{Z}^2} F_X(v)$ with the corresponding sum for the Fourier transform. In view of (6.4–6.7), we conclude

(6.8) $$\sum_{n \geq 0} r_2(n) f_X(\sqrt{n}) = \pi X + \sqrt{X} \sum_{n \geq 1} n^{-1/2} r_2(n) \psi(\delta\sqrt{n}) J_1(2\pi\sqrt{nX}).$$

Recall that the difference between $f_X$ and the characteristic function of the interval $[0, \sqrt{X}]$ is supported on the interval $[\sqrt{X} - \delta, \sqrt{X} + \delta]$ and has absolute value $\leq 1$. Thus

(6.9) $$\Delta(X) \leq B_1(X) + B_2(X) + B_3(X),$$

where

(6.10) $$\begin{aligned} B_1(X) &= \sum \{ r_2(n) \mid \sqrt{X} - \delta \leq \sqrt{n} \leq \sqrt{X} + \delta \}, \\ B_2(X) &= \sqrt{X} \sum_{1 \leq n \leq \delta^{-2}} \frac{r_2(n)}{\sqrt{n}} \left| \psi(\delta\sqrt{n}) J_1(2\pi\sqrt{nX}) \right|, \\ B_3(X) &= \sqrt{X} \sum_{\delta^{-2} < n < \infty} \frac{r_2(n)}{\sqrt{n}} \left| \psi(\delta\sqrt{n}) J_1(2\pi\sqrt{nX}) \right|. \end{aligned}$$

To establish (6.2), it suffices to bound each of these three terms by a multiple of $X^{1/3+\epsilon}$, for every $\epsilon > 0$.

In the case of $B_1(X)$, this follows – with room to spare – from the definition $\delta = X^{-1/6}$ and the bound $r_2(n) = O(n^\epsilon)$; cf. (3.7). For $B_2(X)$, we also use the stationary phase estimate

(6.11) $$|J_1(x)| = O(x^{-1/2})$$

[6, p. 972] and the boundedness of $\psi$:

(6.12) $$B_2(X) = O\left(X^{1/4} \sum_{1 \leq n \leq \delta^{-2}} n^{-3/4+\epsilon}\right) = O(X^{1/3+\epsilon}).$$

We argue similarly in the case of $B_3(X)$, but now use the rapid decay of $\psi$. In particular $\psi(\delta\sqrt{n}) = O(n^{1/2}X^{-1/6})^{-1/2-3\epsilon}$, hence

$$(6.13) \qquad B_3(X) = O\Big(X^{1/3+\epsilon/2}\sum_{n>\delta^{-2}} n^{-1-\epsilon/2}\Big) = O(X^{1/3+\epsilon}).$$

That completes the argument.

With more effort, one can remove the $\epsilon$ from these bounds and get Voronoi's result (6.1). A more serious limitation of the argument above is the use of (6.11) to estimate the size of $J_1$. Various authors have done better by taking advantage of the oscillation of $J_1$ [7, 8, 13, 16, 17]; the current record bound for $\Delta(X)$ is $O(X^{23/73+\epsilon})$ [10]. In the other direction, Hardy [7] used the oscillation to show that $\Delta(X)$ is as big as a constant times $X^{1/4}$ infinitely often. That makes $O(X^{1/4+\epsilon})$ the best possible estimate for the error term, conjectured by Hardy and Landau [8] in 1924. This is analogous to the Riemann hypothesis and the prime number theorem, where the optimal bound on the error term $\pi(X) - \text{Li}(X)$ is thought to be $O(|X|^{1/2}\log X)$, roughly the square root of the trivial bound $O(X/\log X)$. However, the issues are different: the circle problem is one in analysis, in that it asks for a sharp bound for a sum when a harsh cutoff function is introduced; formula (6.8) already gives precise estimates for smoothed sums. In the case of the Riemann hypothesis, the optimal bounds are not known even for smoothed sums.

Stephen D. Miller  
Department of Mathematics  
Hill Center-Busch Campus  
Rutgers University  
110 Frelinghuysen Rd  
Piscataway, NJ 08854-8019  
miller@math.rutgers.edu

Wilfried Schmid  
Department of Mathematics  
Harvard University  
Cambridge, MA 02138  
schmid@math.harvard.edu